\documentclass[12pt,a4paper]{amsart}

\usepackage{enumerate}
\usepackage{xcolor}
\usepackage{amsmath,amssymb}
\usepackage{hyperref}
\usepackage[hmargin=3cm,tmargin=4.5cm,bmargin=3cm]{geometry}
\usepackage{pstricks,pst-node,pst-plot,multido}

\renewcommand{\baselinestretch}{1.2}

\newtheorem*{MT}{Main Theorem}{\bf}{\it}
\newenvironment{Proof}[1][Proof]{\begin{trivlist}
\item[\hskip \labelsep {\itshape #1}.]}{\end{trivlist}}

\newtheorem{lemma}{Lemma}
\newtheorem{proposition}{Proposition}
\theoremstyle{remark}
\newtheorem{remark}{Remark}
\theoremstyle{definition}
\newtheorem{definition}{Definition}

\DeclareMathOperator{\Supp}{Supp} \DeclareMathOperator{\Bk}{Bk}

 \DeclareMathOperator{\NS}{NS}

\DeclareMathOperator{\NE}{NE}

\begin{document}

\title{Log del Pezzo Surfaces of Rank 2 and Cartier Index 3 with a Unique Singularity}

\author{Fei Wang}

\address{
  Fei Wang\\
  Department of Mathematics\\
  National University of Singapore\\
  10 Lower Kent Ridge Road\\
  Singapore 119076
}

\email{matwf@nus.edu.sg}



\maketitle

\begin{abstract}Log del Pezzo surfaces play the role of the opposite of
surfaces of general type. We will completely classify all the log
del Pezzo surfaces of rank 2 and Cartier
index 3 with a unique singularity.
\end{abstract}

The open log del Pezzo surfaces of rank one are discussed by
Miyanishi and Tsunoda in \cite{miyanishi2001open},
\cite{miyanishi1984logarithmic}, \cite{miyanishi1984noncomplete};
and the (complete) log del Pezzo surfaces of rank one are studied by
Kojima \cite{kojima1999logarithmic}, \cite{kojima2003rank}, Zhang
\cite{zhang1988logarithmic}, \cite{zhang1989logarithmic}. Alexeev
and Nikulin give the classification of the log del Pezzo surfaces of
index $\leq 2$ in \cite{alekseev1988classification}, and Nakayama
gives a geometrical classification without using the theory of K3
lattices in \cite{nakayama2007classification}.

\begin{definition}[{\cite[Definition~1]{zhang1989logarithmic}}]\label{d1.2}
  \rm Let $\bar X$ be a normal projective surface with only quotient singularities.
  Then $\bar X$ is called a \emph{logarithmic} (abbr.\ \emph{log})
  \emph{del Pezzo surface} if its anti-canonical divisor $-K_{\bar X}$ is an
  ample $\mathbb Q$-Cartier divisor.
\end{definition}

The smallest positive integer $I$ such that $IK_{\bar X}$ is a
Cartier divisor is called the \emph{Cartier index} of $\bar X$, and
the Picard number $\rho(\bar X)$ is called the \emph{rank} of $\bar
X$. For notations and terminologies, we refer to
Section~\ref{sec1.2}. In the present article, we will give the
complete classification of the log del Pezzo surfaces of rank 2 and
Cartier index 3 with a unique singularity by proving the following.

\begin{MT}\label{theorem1}
  Let $\bar X$ be a log del Pezzo surface with a unique singularity $x_0$,
  and $(X,D)$ the minimal resolution. Suppose that $\bar X$ has rank~$2$ and Cartier
  index~$3$. Then the following assertions hold:

  {\rm 1)} There is a contraction $\pi: \bar X\to \bar Y$ of an irreducible curve $\bar C$ on $\bar X$ to a log del Pezzo surface $\bar Y$ of rank $1$.
  Let $C$ be the proper transform of $\bar C$ on $X$. Then $C$ is a $(-1)$-curve.

  {\rm 2)} The weighted dual graph of $C+D$ is of one of the $29$
  configurations in Figure~\ref{fig1.6}. Moreover, they are all realizable.
\end{MT}

\section{Preliminaries}\label{sec1.2}

We work on an algebraically closed field of characteristic zero.

\begin{definition}[{\cite[Definition~0.2.10]{kawamata1987introduction}\label{d1.3}}]
  \rm Let $\bar X$ be a normal variety. Then $\bar X$ is said to have \emph{log terminal singularities} if

  {\rm 1)} the canonical divisor $K_{\bar X}$ is a $\mathbb Q$-Cartier divisor, i.e., $mK_{\bar X}$ is a Cartier divisor for some $m\in \mathbb Z^+$, and

  {\rm 2)} there exists a resolution of singularities $f: X\to \bar X$ with irreducible exceptional divisors $\{D_j\}_{j=1}^n$ such that $D:=\sum_{j=1}^n D_j$ is a simple normal crossing divisor, and that
        \begin{equation*}
          K_X= f^*(K_{\bar X})+ \sum_{j=1}^n \alpha_j D_j
        \end{equation*}
        for some $\alpha_j\in \mathbb Q$ with $\alpha_j>-1$.
\end{definition}

\medskip

\begin{lemma}[{\cite[Theorem~9.6]{kawamata1988crepant}, \cite[\S4.1]{miyanishi2001open}}]\label{l1.4}
  Suppose $\bar X$ is a normal surface. Then $\bar X$ has only log terminal singularities if and only if $\bar X$ has only quotient singularities. Moreover, if this is the case, let $ X\to \bar X$ be the minimal resolution, then each irreducible exceptional curve is a nonsingular rational curve.
\end{lemma}

Recall that a del Pezzo surface is a normal surface with ample
anti-canonical divisor. It follows from Definition~\ref{d1.3} and
Lemma~\ref{l1.4} that, the log del Pezzo surface as in
Definition~\ref{d1.2} is equivalent to ``the del Pezzo surface with
only log terminal singularities''.

\begin{remark}\label{r1.5}
  Let $\bar X$ be a log del Pezzo surface. Since $\dim \bar X=2$, in Definition~\ref{d1.3} we can take $f: X\to \bar X$ to be the minimal solution. Then $\alpha_j\leq 0$ for all $j$. It follows that $D^{\#}:=-\sum_{j=1}^n \alpha_jD_j$ is an effective $\mathbb Q$-Cartier divisor, and
  $f^*(K_{\bar X})= K_X+D^{\#}$. If $\alpha_k=0$ for some $k$, then $\alpha_j=0$ for all $D_j$ in the connected component of $D$ containing $D_k$
  (\cite[Proposition~4-6-2]{matsuki2002introduction}).
  If $D^{\#}=0$, then $f^*(K_{\bar X})=K_X$ and $\bar X$ is a Gorenstein log del Pezzo surface, which is completely classified in \cite{ye2002gorenstein}. The case when $\bar X$ has index $2$ is classified
  in \cite{alekseev1988classification} and \cite{nakayama2007classification}.
\end{remark}

\medskip

\begin{lemma}[{cf.\ \cite[Lemma~1.1]{zhang1989logarithmic}}]\label{l1.6}
  Let $\bar X$ be a log del Pezzo surface. With the notations in Remark~\ref{r1.5}, we have the following assertions:

  {\rm 1)} $-(K_X+D^{\#})\cdot C\geq 0$ for every irreducible curve $C$ on $X$, and the equality holds if and only if $C\subseteq \Supp(D)$.

  {\rm 2)} If $C\nsubseteq \Supp(D)$ is an irreducible curve on $X$ with negative self-intersection number, then $C$ is a $(-1)$-curve.

  {\rm 3)} $\rho(X)=n+\rho(\bar X)$.
\end{lemma}

\begin{Proof}
  1) Note that $f$ is birational. Since $-K_{\bar X}$ is ample,
  \begin{equation*}
    -(K_X+D^{\#})\cdot C= -f^*(K_{\bar X})\cdot C= -K_{\bar X}\cdot f_*(C)\geq 0.
  \end{equation*}
  The equality holds if and only if $f_*(C)$ is a point, i.e., $C\subseteq \Supp(D)$.

  2) Suppose $C\nsubseteq \Supp(D)$. Then by (1) and the adjunction formula,
  \begin{equation*}
    0<-(K_X+D^{\#})\cdot C \leq -K_X\cdot C = 2+C^2-2p_a(C)\leq 2+C^2\leq 1.
  \end{equation*}
  It follows that $C^2=-1$ and $p_a(C)=0$. So $C$ is a $(-1)$-curve.

  3) $\NS_{\mathbb Q}(X):=\NS(X)\otimes_{\mathbb Z}\mathbb Q$ is generated by $f^*(\NS_{\mathbb Q}(\bar X))$ and
  $\{D_j\}_{j=1}^n$. \qed
\end{Proof}

\medskip

In \cite{kojima1999logarithmic}, $(X,D)$ is assumed to be
\emph{almost minimal}, and we will show in the following that the
minimal resolution of every log del Pezzo surface of rank $1$ is
almost minimal. Hence, we can use the classification in the paper
for our discussion in Sections~\ref{sec1.3}--\ref{sec1.5}.

\medskip

\begin{definition}[{\cite[\S3.11]{miyanishi2001open}\label{d1.7}}]
  \rm Let $\bar X$ be a surface and $(X,D)\to \bar X$ the minimal resolution. With the notations in Remark~\ref{r1.5}, let $\Bk(D)=D-D^{\#}$. Then $(X,D)$ is called \emph{almost minimal} if for every irreducible curve $C$ on $X$ either

  {\rm 1)} $(K_X+D^{\#})\cdot C\geq 0$; or

  {\rm 2)} the intersection matrix of $C+\Bk(D)$ is not negative definite.
\end{definition}

\begin{lemma}\label{l1.8}
  Let $\bar X$ be a log del Pezzo surface of rank $1$. Then its minimal resolution $(X,D)$ is almost minimal.
\end{lemma}

\begin{Proof}
  Suppose there exists an irreducible curve $E$ on $X$ such that $E\cdot (K_X+D^{\#})<0$
  and the intersection matrix of $E+\Bk(D)$, i.e., of $E+D$, is negative definite.

  Let $\bar E=f_*(E)$. Since $0>E\cdot f^*(K_{\bar X})=\bar E\cdot K_{\bar X}$, $\bar E$ is a curve on $\bar X$. Recall that $\rho(\bar X)=1$. We
  can write $\bar E\equiv rK_{\bar X}$ for some $r\in \mathbb Q$. Then $(\bar E)^2=r^2(K_{\bar X})^2\geq 0$.

  On the other hand,
  \begin{equation*}
    f^*(\bar E)= E+ \sum_{j=1}^n \beta_j D_j
  \end{equation*}
  for some $\beta_j\in \mathbb Q$. Let $H= \sum_{j=1}^n \beta_j D_j$. Then
  \begin{equation*}
    (\bar E)^2= (f^*(\bar E))^2=(E+H)^2<0,
  \end{equation*}
  because the intersection matrix of $E+D$ is negative definite. This leads to a
  contradiction. \qed
\end{Proof}

\section[Types of Weighted Dual Graphs]{The Types of Weighted Dual Graphs of $D$}\label{sec1.3}

In this section, we assume that $\bar X$ is a log del Pezzo surface
of Cartier index 3 with a unique singularity $x_0$, and use the
notations in Section~\ref{sec1.2}. Note that the dual graph of the
exceptional divisor $D$ is of A-D-E Dynkin's type. We are going to
determine all the possible types of the weighted dual graphs of $D$.
\begin{figure}[h!]
\begin{center}\footnotesize \psset{dotstyle=o,linewidth=0.5pt}
  \begin{tabular}{|c|c|c|}
    \hline
    No. & Weighted Dual graph of $D$ & Size \\
    \hline
    I &
    \begin{pspicture*}(-1,-0.5)(1,0.3)
      \psdot(0,0)
      \uput[-90](0,0){$-3$}
    \end{pspicture*}
    & $n=1$ \\
    \hline
    II &
    \begin{pspicture*}(-1,-0.5)(1,0.3)
      \psdot(0,0)
      \uput[-90](0,0){$-6$}
    \end{pspicture*}
    & $n=1$ \\
    \hline
    III &
    \begin{pspicture*}(-1,-0.5)(2,0.3)
      \psline[arrows=o-o](0,0)(1,0)
      \uput[-90](0,0){$-2$}
      \uput[-90](1,0){$-5$}
    \end{pspicture*}
    & $n=2$ \\
    \hline
    IV &
    \begin{pspicture*}(-1,-0.5)(3,0.3)
      \psline[arrows=o-](0,0)(1,0)
      \psline[arrows=o-o](1,0)(2,0)
      \uput[-90](0,0){$-2$}
      \uput[-90](1,0){$-4$}
      \uput[-90](2,0){$-2$}
    \end{pspicture*}
    & $n=3$ \\
    \hline
    V &
    \begin{pspicture*}(-1,-0.5)(4,0.3)
      \psline[arrows=o-](0,0)(1,0)
      \psline[linestyle=dashed,arrows=o-](1,0)(2,0)
      \psline[arrows=o-o](2,0)(3,0)
      \uput[-90](0,0){$-4$}
      \uput[-90](1,0){$-2$}
      \uput[-90](2,0){$-2$}
      \uput[-90](3,0){$-4$}
    \end{pspicture*}
    & $2\leq n\leq 10$ \\
    \hline
    VI & \begin{pspicture*}(-1,-0.5)(5,0.3)
      \psline[arrows=o-](0,0)(1,0)
      \psline[arrows=o-](1,0)(2,0)
      \psline[linestyle=dashed,arrows=o-](2,0)(3,0)
      \psline[arrows=o-o](3,0)(4,0)
      \uput[-90](0,0){$-2$}
      \uput[-90](1,0){$-3$}
      \uput[-90](2,0){$-2$}
      \uput[-90](3,0){$-2$}
      \uput[-90](4,0){$-4$}
    \end{pspicture*}
    & $3\leq n\leq 9$ \\
    \hline
    VII & \begin{pspicture*}(-1,-0.5)(6,0.3)
      \psline[arrows=o-](0,0)(1,0)
      \psline[arrows=o-](1,0)(2,0)
      \psline[linestyle=dashed,arrows=o-](2,0)(3,0)
      \psline[arrows=o-](3,0)(4,0)
      \psline[arrows=o-o](4,0)(5,0)
      \uput[-90](0,0){$-2$}
      \uput[-90](1,0){$-3$}
      \uput[-90](2,0){$-2$}
      \uput[-90](3,0){$-2$}
      \uput[-90](4,0){$-3$}
      \uput[-90](5,0){$-2$}
    \end{pspicture*}
    & $4\leq n\leq 9$ \\
    \hline
    VIII & \begin{pspicture*}(-1,-0.5)(4,1.4)
      \psline[arrows=o-](0,0)(1,0)
      \psline[arrows=o-,linestyle=dashed](1,0)(2,0)
      \psline[arrows=-o](2,0)(3,0)
      \psline[arrows=o-o](2,0)(2,1)
      \uput[-90](0,0){$-4$}
      \uput[-90](1,0){$-2$}
      \uput[-90](2,0){$-2$}
      \uput[-90](3,0){$-2$}
      \uput[180](2,1){$-2$}
    \end{pspicture*}
    & $4\leq n\leq 9$ \\
    \hline
    IX & \begin{pspicture*}(-1,-0.5)(5,1.4)
      \psline[arrows=o-](0,0)(1,0)
      \psline[arrows=o-](1,0)(2,0)
      \psline[arrows=o-,linestyle=dashed](2,0)(3,0)
      \psline[arrows=-o](3,0)(4,0)
      \psline[arrows=o-o](3,0)(3,1)
      \uput[-90](0,0){$-2$}
      \uput[-90](1,0){$-3$}
      \uput[-90](2,0){$-2$}
      \uput[-90](3,0){$-2$}
      \uput[-90](4,0){$-2$}
      \uput[180](3,1){$-2$}
    \end{pspicture*}
    & $4\leq n\leq 8$ \\
    \hline
  \end{tabular}
\end{center}\caption{Weighted Dual graph of $D$}\label{fig1.1}
\end{figure}

Let $a_j=-\alpha_j$. Then $f^*(K_{\bar X})= K_X+\sum_{j=1}^n a_jD_j$
for some $0<a_j<1$. It is given that $3K_{\bar X}$ is a Cartier
divisor, so is $\sum_{j=1}^n a_jD_j$. Therefore, $a_j\in
\{1/3,2/3\}$ for all $j$. Note that for each $i=1,\ldots,n$,
\begin{equation*}
  0=f^*(K_{\bar X})\cdot D_i=\left(K_X+ \sum_{j=1}^n a_jD_j\right)\cdot D_i= -2-(D_i)^2+ \sum_{j=1}^n a_j(D_i\cdot D_j).
\end{equation*}
That is,
\begin{equation*}
  \sum_{j=1}^n a_j(D_i\cdot D_j)= 2+(D_i)^2,\quad i=1,\ldots,n.
\end{equation*}
Using these results, we can show that

\begin{proposition}\label{p1.9}
  Let $\bar X$ be a log del Pezzo surface of Cartier index $3$ with a unique singularity, and $(X,D)$
  its minimal resolution. Then

  {\rm 1)} the weighted dual graph of $D$ is of one of the nine cases listed
  in the second column of Figure~\textup{\ref{fig1.1}}, and

  {\rm 2)} the possible sizes of $D$ are given in the third column
  of Figure~\textup{\ref{fig1.1}}.
\end{proposition}

We will leave the proof of (2) in Section~\ref{sec1.4}.

\begin{Proof}[Proof of Proposition~\textup{\ref{p1.9} (1)}]
Consider the two cases:

\paragraph{\bf Type A}
Suppose that $D$ is a linear chain $D_1-D_2-\cdots-D_n$.

If $n=1$, then $a_1(D_1)^2=2+(D_1)^2$. When $a_1=1/3$, $(D_1)^2=-3$,
and $D$ is given by I of Figure~\ref{fig1.1}; when $a_1=2/3$,
$(D_1)^2=-6$, and $D$ is given by II.

Suppose $n\geq 2$. Then for all $i=2,\ldots,n$,
$a_{i-1}+a_i(D_i)^2+a_{i+1}=2+(D_i)^2$. This implies
$2-a_{i-1}-a_{i+1}=(D_i)^2(a_i-1)\geq -2(a_i-1)$, i.e.,
\begin{equation*}
  a_i\geq \dfrac{1}{2}(a_{i-1}+a_{i+1}).
\end{equation*}
Moreover, the equality holds if and only if $(D_i)^2=-2$.

If $a_i=1/3$ for some $i=2,\ldots,n-1$, then $a_{i-1}+a_{i+1}\leq
2/3$ and thus $a_{i-1}=a_{i+1}=1/3$; consequently $a_j=1/3$ for all
$j=1,\ldots,n$. In particular, $1/3\,(D_1)^2+1/3=2+(D_1)^2$.
However, this would imply that $(D_1)^2=-5/2\notin \mathbb Z$, a
contradiction. So $a_i=2/3$ for some $i=2,\ldots,n-1$. If $i\leq
n-2$,  then $a_{i+1}\geq \frac{1}{2}(a_i+a_{i+2})\geq
\frac{1}{2}(\frac{2}{3}+\frac{1}{3})=1/2$, and then $a_{i+1}=2/3$.
It follows by induction that $a_j=2/3$ for all $j=i,\ldots,n-1$; and
similarly $a_j=2/3$ for all $j=2,\ldots,i$. We consider three cases:

(i) $a_j=2/3$ for all $j=1,\ldots,n$. Then $(D_1)^2=(D_n)^2=-4$  and
$(D_j)^2=-2$ for $j=2,\ldots,n-1$. This is given by V of
Figure~\ref{fig1.1}.

(ii) $a_1=1/3$ and $a_j=2/3$ for all $j=2,\ldots,n$. For this case,
if $n=2$, then $(D_1)^2=-2$ and $(D_2)^2=-5$, which is given by III;
if $n\geq 3$, then $(D_2)^2=-3$, $(D_n)^2=-4$ and $(D_j)^2=-2$ for
all other $j$, which is given by VI of Figure~\ref{fig1.1}.

(iii) $a_1=a_n=1/3$ and $a_j=2/3$ for all $j=2,\ldots,n-1$.  It is
impossible if $n=2$. If $n=3$, then $(D_1)^2=(D_3)^2=-2$ and
$(D_2)^2=-4$, which is given by IV; if $n\geq 4$, then
$(D_2)^2=(D_{n-1})^2=-3$ and $(D_j)^2=-2$ for all other $j$, which
is given by VII.

\paragraph{\bf Type D and E} Suppose that $D$ is a fork.
Let $D_3$ be the center of the fork. It intersects with three
components, say $D_1,D_2$ and $D_4$. Then
$a_1+a_3+a_4+a_2(D_2)^2=2+(D_2)^2$. There are two cases:

(i) If $(D_3)^2\leq -3$, then $1\geq
2-a_1-a_2-a_4=(D_3)^2(a_3-1)\geq (-3)(1/3)=1$.  We have
$a_1=a_2=a_4=1/3$, $a_3=2/3$ and $(D_3)^2=-3$. If $D_4$ intersects
with, say, $D_5$, then $2/3+a_5+1/3\,(D_4)^2=2+(D_4)^2$ implies
$(D_4)^2=(3/2)a_5-2\geq -3/2$, a contradiction. So $D_4$ is the end
of a twig, and the same is true for $D_1$ and $D_2$. Therefore, for
this case $n=4$ and $(D_1)^2=(D_2)^2=(D_4)^2=-2$. The weighted dual
graph is by IX $(n=4)$.

(ii) If $(D_3)^2=-2$, then $a_1+a_2+a_4=2a_3$. It follows that
$a_3=2/3$ and $a_1+a_2+a_4=4/3$. After the relabeling if necessary,
we have $a_1=a_2=1/3$ and $a_4=2/3$. Using the same argument as
above, $D_1$ and $D_2$ are twigs of $D$ consisting of a single
$(-2)$-curve. We are left to determine the last twig of $D$:
$\displaystyle{{D_1}\atop {D_2}}\!\!>\!D_3-D_4-\cdots -D_n$. Using
the same argument as in the case of linear chain, it follows by
induction that $a_j=2/3$ for all $j=4,\ldots,n-1$. There are two
cases:

(ii.a) $a_1=a_2=1/3$ and $a_j=2/3$ for all $j=3,4,\ldots,n$.  Then
$(D_n)^2=-4$ and $(D_j)^2=-2$ for all $j=1,\ldots,n-1$. This is
given by VIII of Figure~\ref{fig1.1}.

(ii.b) $a_1=a_2=a_n=1/3$ and $a_j=2/3$ for all $j=3,4,\ldots,n-1$.
Then $n\geq 5$, $(D_{n-1})^2=-3$ and $(D_j)^2=-2$ for all $j\neq
n-1$. This is given by IX $(n\geq 5)$. \qed
\end{Proof}

\section{Contraction}\label{sec1.4}

From now on, we assume that $\bar X$ is a log del Pezzo surface of
rank 2 and Cartier index 3 with a unique singularity $x_0$. Since
$K_{\bar X}$ is not numerically effective, by cone theorem, there is
a $K_{\bar X}$-negative extremal ray  $R\subseteq
\overline{\NE}(\bar X)$. Let $\pi: \bar X\to \bar Y$ be the
contraction of $R$. Then $\bar Y$ is a normal projective variety of
$\dim \bar Y\leq 2$ and $\pi$ has connected fibers. We will consider
the three possibilities according to the dimension of $\bar Y$.

\medskip

\emph{Case 1}: $\dim \bar Y=0$. It follows that $N_1(\bar X)$ is
generated by some $[\bar C]\in R$, and thus $\rho(\bar X)=1$. But we
assumed that $\rho(\bar X)=2$, a contradiction.

\medskip

\emph{Case 2}: $\dim \bar Y=1$. Then $\bar Y$ is a nonsingular
curve. By \cite[Lemma~1.3]{gurjar1994pi1}, every log del Pezzo
surface is a rational surface. Then it follows from L\"{u}roth's
theorem that the base $\bar Y$ is rational. Therefore, $Y\cong
\mathbb P^1$. We claim that

\begin{lemma}\label{l1.10}
  With the notations above, every fiber of the contraction $\pi: \bar X\to \bar Y$ is irreducible.
\end{lemma}
\begin{Proof}
  Since $\bar Y$ is nonsingular, the contraction $\pi: \bar X\to \bar Y$ is flat, and thus every fiber has pure dimension $1$.
  For any point $y\in \bar Y$, let $\bar F=\pi^{-1}(y)$. Suppose $\bar
F$ is reducible. Since $\bar F$ is connected, we may choose
irreducible components $\bar F_1$ and $\bar F_2$ of $\bar F$ such
that $\bar F_1\cdot \bar F_2\geq 1$. On the other hand, $\bar
F_1\equiv a\bar F_2\in R$ for some $a>0$. Then by Zariski's lemma
\cite[Lemma~8.2]{barth2004compact}, $\bar F_1\cdot \bar F_2= a(\bar
F_2)^2<0$, a contradiction. \qed
\end{Proof}

We continue the discussion of $\dim \bar Y=1$. Let $y_0=\pi(x_0)$
and $\bar C=\pi^{-1}(y_0)$. Then $x_0\in \bar C$, and by Zariski's
lemma, $(\bar C)^2=0$. Take $f: (X,D)\to \bar X$ to be the minimal
resolution, and $C$ the proper transform of $\bar C$ with respect to
$f$. Then $C+D=(\pi\circ f)^{-1}(y_0)$. By Zariski's lemma again,
$C^2<0$, and thus $C$ is a $(-1)$-curve by Lemma~\ref{l1.6}.

Let $y\in \bar Y\backslash \{y_0\}$, $\bar F:=\pi^{-1}(y)$ and $F$
the proper transform of $\bar F$ with respect to $f$. Then
$F=(\pi\circ f)^{-1}(y)$. So $F^2=0$ and $F\cdot D^{\#}=0$. We have
\begin{equation*}
  0>\bar F\cdot K_{\bar X}= F\cdot (K_X+D^{\#})=F\cdot K_X.
\end{equation*}
Then by adjunction formula, $2p_a(F)-2=F\cdot (F+K_X)=F\cdot K_X<0$,
and thus $p_a(F)=0$. By Lemma~\ref{l1.10}, $F$ is irreducible; so
$F\cong \mathbb P^1$.

Let $F_0$ be the singular fiber of the the $\mathbb P^1$-fibration
$\pi\circ f: X\to \bar Y$ over $y_0$. Then $\Supp(F_0)=C+D$. After
contracting $C$ and consecutively $(-1)$-curves in $C+D$, $C+D$
becomes $\mathbb P^1$. In particular, note that $D$ is connected and
$C+D$ is a connected simple normal crossing divisor, we have $C\cdot
D=1$. Moreover,
\begin{equation}\label{eqn1.1}
  2+n=\rho(X)=10-(K_X)^2.
\end{equation}

\medskip

\emph{Case 3}: $\dim \bar Y=2$. Then $\pi: \bar X\to \bar Y$ is
birational and the exceptional curve is irreducible
\cite[Proposition~2.5]{kollar1998birational}, denoted by $\bar C$.
Let $C$ be the proper transform of $\bar C$ with respect to the
minimal resolution $f: (X,D)\to \bar X$.

Note that $\pi\circ f: X\to \bar Y$ contracts $C$ into a point.  By
negative definiteness theorem, $C^2<0$. So by Lemma~\ref{l1.6}, $C$
is a $(-1)$-curve. By
\cite[Proposition~5-1-6]{kawamata1987introduction}, $\bar Y$ is
$\mathbb Q$-factorial, and it is either smooth or it has a unique
log terminal singularity $y_0=\pi(x_0)$. By taking $H=-K_{\bar X}$
in Lemma~\ref{l1.11} below, $-K_{\bar Y}$ is ample. Therefore, $\bar
Y$ is either a smooth del Pezzo surface or a log del Pezzo surface
with a unique singularity $y_0$. Recall that $\rho(\bar Y)=1$. If
$\bar Y$ is smooth, then $\bar Y\cong \mathbb P^2$, the projective
plane.

\begin{lemma}\label{l1.11}
  With the notations as above, for any ample divisor $H$ on $\bar X$, $\pi_*(H)$ is ample.
\end{lemma}

\begin{Proof}
  Let $\bar H=\pi_*(H)$. Then by projection formula $H=\pi^*(\bar H)+a\bar C$ for some $a\in \mathbb R$. Suppose $x_0\in \bar C$. Since $f^{-1}(\bar C)= \Supp(C+D)$ and that
the
   intersection matrix of $C+D$ is negative definite, $(\bar C)^2= (f^*(\bar C))^2<0$.
   If $x_0\notin \bar C$, then $(\bar C)^2= C^2=-1$. For either case,
  \begin{equation*}
    0<H^2=(\pi^*(\bar H)+a\bar C)^2= (\pi^*(\bar H))^2+ a^2(\bar C)^2\leq (\pi^*(\bar H))^2= (\bar H)^2.
  \end{equation*}
Let $\bar E$ be an irreducible curve on $\bar Y$ and
  $\bar E'$ the proper transform of $\bar E$ with respect to $\pi$. Then $\pi^*(\bar E)=\bar E'+b\bar C$ for some $b\in \mathbb R$. We can compute that
  \begin{equation*}
    0=\bar C\cdot \pi^*(\bar E)= \bar C\cdot \bar E'+b(\bar C)^2 \geq b(\bar C)^2.
  \end{equation*}
  So $b\geq 0$. Then
  \begin{equation*}
    \bar H\cdot \bar E = H\cdot \pi^*(\bar E)= H\cdot (\bar E'+b\bar
    C)= H\cdot \bar E'+b(H\cdot \bar C)\geq H\cdot \bar E'>0.
  \end{equation*}
  By Nakai-Moishezon criterion, $\bar H$ is an ample divisor on $\bar
  Y$. \qed
\end{Proof}

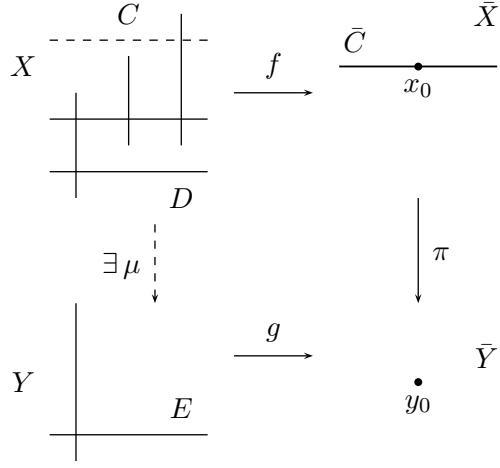
\begin{figure}[h!]
\begin{center}\psset{unit=0.7,linewidth=0.5pt}\everypsbox{\small}
\begin{pspicture*}(-5,-5)(5.5,4.5)
  \psline(-3.5,1)(-0.5,1)
  \psline(-3,0.5)(-3,2.5)
  \psline(-3.5,2)(-0.5,2)
  \psline(-2,1.5)(-2,3.2)
  \psline(-1,1.5)(-1,4)
  \psline[linestyle=dashed, dash=3pt 3pt](-3.5,3.5)(-0.5,3.5)

  \psline[linewidth=0.7pt](2,3)(5,3)
  \psdots(3.5,3)
  \uput[-90](3.5,3){$x_0$}

  \psline(-3.5,-4)(-0.5,-4)
  \psline(-3,-4.5)(-3,-1.5)

  \psdots(3.5,-3)
  \uput[-90](3.5,-3){$y_0$}

  \pcline[arrows=->](0,2.5)(1.5,2.5)
  \taput{$f$}
  \pcline[arrows=->](0,-2.5)(1.5,-2.5)
  \taput{$g$}
  \pcline[arrows=->](3.5,0.5)(3.5,-1.5)
  \trput{$\pi$}
  \pcline[arrows=->, linestyle=dashed,dash=3pt 3pt](-1.5,0)(-1.5,-1.5)
  \tlput{$\exists\, \mu$}

  \rput(-4,3){$X$}
  \rput(4.8,4){$\bar X$}
  \rput(-4,-3){$Y$}
  \rput(4.8,-2.5){$\bar Y$}
  \rput(-2,4){$C$}
  \rput(2.3,3.5){$\bar C$}
  \rput(-1,0.5){$D$}
  \rput(-1,-3.5){$E$}
\end{pspicture*}\caption{Divisorial Contraction}\label{fig1.2}
\end{center}
\end{figure}

We continue the discussion of $\dim \bar Y=2$. Let $g: Y\to \bar Y$ be the minimal resolution. Then $\pi\circ f$
factors through $Y$; that is, there is a proper birational morphism
$\mu: X\to Y$ such that $g\circ \mu=\pi\circ f$ as illustrated in
Figure~\ref{fig1.2}. We see that $\mu: X\to Y$ is the composite of
blow-downs of $(-1)$-curves. More precisely, it is the contraction
of $C$ and consecutive $(-1)$-curves in $C+D$.

Let $y_0=f(x_0)$. If $\bar Y\cong \mathbb P^2$, then $Y=\bar Y$ and
$\mu(C+D)=y_0$. Suppose $\bar Y$ is a log del Pezzo surface of rank
$1$ with a unique singularity $y_0$. Then $Y$ can be further
contracted along $(-1)$-curves into the Hirzebruch surface $\mathbb
F_r$ for some $r\geq 0$ \cite[Theorem~2.1, 3.1,
4.1]{kojima1999logarithmic}. For either case,
\begin{equation}\label{eqn1.2}
  2+n=\rho(X)=10-(K_X)^2.
\end{equation}

We can now determine the size of the weighted dual graphs of $D$ in
Figure~\ref{fig1.1}.

\begin{Proof}[Proof of Proposition~\textup{\ref{p1.9} (2)}]
Recall that $-K_{\bar X}$ is ample. In particular,
\begin{align*}
  0  <(K_{\bar X})^2&  = \pi^*(K_{\bar X})\cdot \pi^*(K_{\bar X}) = K_X\cdot \pi^*(K_{\bar X}) \\
  & = K_X\cdot \left(K_X+ \sum_{j=1}^n a_jD_j\right) \\
  & = (K_X)^2+ \sum_{j=1}^n a_j(-2-(D_j)^2).
\end{align*}
For both the fiber contraction (\ref{eqn1.1}) and the divisorial
contraction (\ref{eqn1.2}),
\begin{equation*}
  2+n=\rho(X)= 10-(K_X)^2<10+\sum_{j=1}^n a_j(-2-(D_j)^2).
\end{equation*}
That is, $n<8+\sum_{j=1}^n a_j(-2-(D_j)^2)$. Recall that
$D^{\#}=\sum_{j=1}^n a_jD_j$ is evaluated explicitly in the proof of
part (1), we can easily compute the possible size $n$ of $D$:
\begin{enumerate}[\quad I.]\itemsep=1mm
  \setcounter{enumi}{4}
  \item $n<8+2/3\cdot 2+2/3\cdot 2\Leftrightarrow n\leq 10$;
  \item $n<8+2/3\cdot 1+2/3\cdot 2\Leftrightarrow n\leq 9$;
  \item $n<8+2/3\cdot 1+2/3\cdot 1\Leftrightarrow n\leq 9$;
  \item $n<8+2/3\cdot 2\Leftrightarrow n\leq 9$;
  \item $n<8+2/3\cdot 1\Leftrightarrow n\leq 8$.
\end{enumerate}
This completes the proof of Proposition~\ref{p1.9} (2). \qed
\end{Proof}

\begin{Proof}[Proof of Main Theorem] 1)
  Suppose $\dim \bar Y=1$. We have seen that $C+D$ can be smoothly contracted to $F\cong \mathbb P^1$ with $F^2=0$
  along $C$ and consecutive $(-1)$-curves in $C+D$. However,
  by verifying all the weighted dual graphs in Figure~\ref{fig1.1}, we see that
  none of them with any $(-1)$-curve can be contracted to such a curve, a contradiction. Therefore, $\dim \bar Y=2$ and $\bar Y$ is a log del Pezzo surface
of rank 1. In particular, as proved in Section~\ref{sec1.4}, $C$ is
a $(-1)$-curve.

\medskip

2) \emph{Case 1.} If $\bar Y$ is smooth, then $Y=\bar Y\cong \mathbb
P^2$ and $C+D$ is contracted to the smooth point $y_0$ along $C$ and
consecutive $(-1)$-curves in $C+D$. In particular, by noting that
$D$ is a simple normal crossing divisor, we have $C\cdot D=1$.

\emph{Case 2.} Suppose $\bar Y$ is not smooth. Then $\bar Y$ is a
log del Pezzo surface with a unique singularity $y_0$. Let $E$ be
the exceptional divisor of the minimal resolution $g:Y\to \bar Y$.
The configuration of $E$ is completely classified in
\cite[Theorem~2.1]{kojima1999logarithmic}.  Recall that the possible
weighted dual graphs of $D$ have been listed in Figure~\ref{fig1.1}.

(i) If $x_0\notin \bar C$, then $C$ is disjoint from $D$, and the
weighted dual graphs of $D$ is the same as that of $E$.

(ii) If $x_0\in \bar C$, then $C+D$ is a connected simple normal
crossing divisor since $E$ is of A-D-E Dynkin's type. Note that $D$
is connected. Then $C\cdot D=1$ and $X\backslash (C\cup D)\cong
Y\backslash E$. We only need to check how $C+D$ is contracted to $E$
along $C$ and consecutive $(-1)$-curves in $C+D$.

By checking all the possible weighted dual graphs of $D$ in
Figure~\ref{fig1.1} and all the possible places of $C$, there are 3
configurations of $C+D$ (VI $(n=5)$ (b), VI $(n=6)$ (b), IX $(n=5)$
(b)) for the case when $\bar Y$ is smooth, and 26 configurations of
$C+D$ for the case when $\bar Y$ is not smooth. They are given in
Figure~\ref{fig1.6}.

According to the discussions above, each of these 29 possible
configurations of $C+D$ can be contracted to $E$ (resp.\ a smooth
point) along $C$ and consecutive $(-1)$-curves in $C+D$. There
exists a log del Pezzo surface $\bar Y$ of rank 1 with a unique
singularity (resp.\ $\bar Y\cong \mathbb P^2$), such that $E$ is the
exceptional divisor of its minimal resolution $Y\to \bar Y$ (resp.\
$Y=\bar Y$). We can construct the surface $X$ by blowing up points
from the corresponding surface $Y$. Let $X\to \bar X$ be the
contraction of $D$. Then $\bar X$ is a projective normal surface of
rank 2 and Cartier index 3 with a unique quotient singularity. We
claim that
\begin{lemma}\label{l1.12}
  For each of the configuration of $C+D$ in Figure~\textup{\ref{fig1.6}}, let $\bar X$ be the surface defined above, then $-K_{\bar X}$ is ample.
\end{lemma}
It follows that $\bar X$ is a log del Pezzo surface of rank 2 and
Cartier index 3 with a unique singularity $x_0$, and $D$ is the
exceptional divisor of its minimal resolution $X\to \bar X$. In
other words, every configuration in Figure~\ref{fig1.6} is
realizable. We have completed the proof of Main Theorem. \qed
\end{Proof}

\section{Ampleness of $-K_{\bar X}$}\label{sec1.5}

In the proof of Main Theorem, for each weighted graph of $C+D$ in
Figure~\ref{fig1.6}, we constructed a normal projective surface
$\bar X$ of rank 2 and Cartier index 3 with a unique quotient
singularity, such that $D$ is the exceptional divisor of its minimal
resolution $X\to \bar X$. In order to prove that $\bar X$ is a log
del Pezzo surface, it remains to show that $-K_{\bar X}$ is ample
(cf.\ Lemma~\ref{l1.12}.)

First of all, we shall evaluate $-K_{\bar X}$. We explore the
notations used in the discussion of the divisorial contraction case
in Section~\ref{sec1.4} (as illustrated in Figure~\ref{fig1.2}).
Recall that $\mu:X\to Y$  is the successive contraction of
$(-1)$-curves in $C+D$. If $\bar Y$ is smooth, then $Y=\bar Y\cong
\mathbb P^2$, and $\mu$ factors through $X\to \mathbb F_1\to Y$. If
$\bar Y$ has a unique singularity, then $Y$ can be further
contracted to the Hirzebruch surface $\mathbb F_r$ for some $r\geq
0$ along $(-1)$-curves \cite[Theorem~3.1,
4.1]{kojima1999logarithmic}.

We can verify the list of configurations in Figure~\ref{fig1.6} to
conclude that

\begin{lemma}\label{l1.13}
  Let $\bar X$ be a log del Pezzo surface of rank $2$ and Cartier index $3$ with a unique singularity, and $(X,D)\to X$ the minimal resolution. Then there exists a $\mathbb P^1$-fibration $X\xrightarrow{\Phi} \mathbb F_r\to \mathbb P^1$ with at most two singular fibers, such that one of the components $D_\ell$ of $D$ is a cross-section, $C$ and the other components of $D$ are contained in the singular fibers.
\end{lemma}

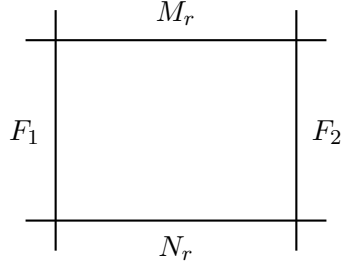
\begin{figure}[h!]
\begin{center}\everypsbox{\small}\psset{unit=0.8}
\begin{pspicture*}(-4,-1.7)(4,2.7)
  \pcline(-2.5,2)(2.5,2)
  \taput{$M_r$}
  \pcline(-2.5,-1)(2.5,-1)
  \tbput{$N_r$}
  \pcline(-2,-1.5)(-2,2.5)
  \tlput{$F_1$}
  \pcline(2,-1.5)(2,2.5)
  \trput{$F_2$}
\end{pspicture*}
\caption{$-K_{\mathbb F_r}$}\label{fig1.3}
\end{center}
\end{figure}

Then $M_r:=\Phi(D_\ell)$ is the minimal section of $\mathbb F_r$. If
there are two singular fibers, let their images in $\mathbb F_r$ be
$F_1$ and $F_2$. If there is only one singular fiber, let its image
in $\mathbb F_r$ be $F_1$ and take $F_2$ to be the image of a
general fiber. Take a section $N_r\sim M_r+rF_1$ which does not
contain the image of any center of blowup. Then $-K_{\mathbb
F_r}=M_r+N_r+F_1+F_2$, which form a circle (Figure~\ref{fig1.3}).

We can decompose $\Phi: X\to \mathbb F_r$ as the composite of
blow-downs $X=X_0\xrightarrow{\phi_1} X_1\to \cdots \to X_{k-1}
\xrightarrow{\phi_k} X_k=\mathbb F_r$. Denote the exceptional curve
of $\phi_i$ by $E_i$, $i=1,\ldots,k$. Then $K_{X_{i-1}}=
\phi_i^*(K_{X_i})+E_i$. Therefore, $-K_X$ can be evaluated explicitly.

Note that $-K_X$ is supported by
$\Delta:=\Phi^{-1}(M_r+N_r+F_1+F_2)$. Let $\Delta_+$ denote the sum
of the irreducible curves which have positive coefficients appearing
in $-K_X$. Note that $\Delta_+$ forms a loop, and every irreducible
curve in $\Delta_+$ has coefficient $1$ appearing in $-K_X$. In
particular, the proper transforms of $M_r,N_r,F_1$ and $F_2$ on $X$
belong to $\Delta_+$.

Recall that in the proof of Proposition~\ref{p1.9} (1), we computed
the unique numbers $a_j\in \{1/3,2/3\}$, $i=1,\ldots,n$, such that
\begin{equation*}
  f^*(K_{\bar X})=K_X+ \sum_{j=1}^n a_jD_j.
\end{equation*}
We can thus evaluate $-f^*(K_{\bar X})$ explicitly.

\begin{figure}[h!]
\begin{center}\psset{dash=3pt 2pt,unit=0.7}\everypsbox{\scriptsize}
\begin{pspicture*}(-0.7,-1.3)(9.8,5.2)
  \psline(-0.4,0)(9.4,0)
  \psline(0,-0.4)(0,1.9)
  \psline(-0.4,1.5)(2.4,1.5)
  \psline(2,1.1)(2,3.4)
  \psline(1.6,3)(3.9,3)
  \psline(3.5,2.6)(3.5,4.9)
  \psline(3.1,4.5)(5.9,4.5)
  \psline(5.5,2.6)(5.5,4.9)
  \psline(5.1,3)(7.4,3)
  \psline(6.7,3.15)(9.3,1.85)
  \psline(9,-0.4)(9,2.4)
  \psline[linestyle=dashed](7,1)(9.8,1)
  \psline[linestyle=dashed](8,1.6)(9,4.5)
  \psline[linestyle=dashed](1,1)(1,4.2)
  \rput(4.5,-0.4){$0$}
  \rput(-0.4,0.75){$-4$}
  \rput(0.4,0.75){$(\frac{1}{3})$}
  \rput(1.5,1.9){$(\frac{1}{3})$}
  \rput(0.6,3.5){$(0)$}
  \rput(2.4,2.3){$(\frac{1}{3})$}
  \rput(2.75,3.4){$(\frac{1}{3})$}
  \rput(3.9,3.75){$(\frac{1}{3})$}
  \rput(4.5,4.9){$(\frac{1}{3})$}
  \rput(5.9,3.75){$(\frac{1}{3})$}
  \rput(6.2,2.6){$(\frac{1}{3})$}
  \rput(7.5,2.2){$(\frac{1}{3})$}
  \rput(8.6,0.5){$(\frac{1}{3})$}
  \rput(9.4,0.5){$-4$}
  \rput(8.3,4){$(0)$}
  \rput(7.5,0.6){$(0)$}
  \rput(6.6,1){$C$}
  \rput(4.5,-1){V $(n=10)$}
\end{pspicture*}
\begin{pspicture*}(-0.7,-1.3)(9.8,5.2)
  \psline(-0.4,0)(9.4,0)
  \psline(0,-0.4)(0,1.9)
  \psline(-0.4,1.5)(2.4,1.5)
  \psline(2,1.1)(2,3.4)
  \psline(1.6,3)(3.9,3)
  \psline(3.5,2.6)(3.5,4.9)
  \psline(3.1,4.5)(6.4,4.5)
  \psline(6,1.6)(6,4.9)
  \psline(5.6,2)(9.4,2)
  \psline(9,-0.4)(9,2.4)
  \psline[linestyle=dashed](7,1)(7,4)
  \psline[linestyle=dashed](8,1)(8,4)
  \psline[linestyle=dashed](1,1)(1,4.2)
  \rput(4.5,-0.4){$0$}
  \rput(-0.4,0.75){$-4$}
  \rput(0.4,0.75){$(\frac{1}{3})$}
  \rput(1.5,1.9){$(\frac{1}{3})$}
  \rput(0.6,3.5){$(0)$}
  \rput(2.4,2.3){$(\frac{1}{3})$}
  \rput(2.75,3.4){$(\frac{1}{3})$}
  \rput(3.9,3.75){$(\frac{1}{3})$}
  \rput(4.75,4.9){$(\frac{1}{3})$}
  \rput(5.6,3.25){$(\frac{1}{3})$}
  \rput(8.5,2.4){$(\frac{1}{3})$}
  \rput(7.4,3.5){$(0)$}
  \rput(8.4,3.5){$(0)$}
  \rput(9.4,1){$(\frac{2}{3})$}
  \rput(8.5,1.6){$-3$}
  \rput(6.6,1){$C$}
  \rput(4.5,-1){VI $(n=9)$}
\end{pspicture*}
\psset{unit=1.15}
\begin{pspicture*}(0.8,-1.4)(9.6,4.8)
  \psline(5.6,0)(9.4,0)
  \psline(9,-0.4)(9,4.4)
  \psline(6.6,4)(9.4,4)
  \psline[linestyle=dashed](6,-0.4)(6,2.4)
  \psline(5.8,1.6)(7.2,4.4)
  \psline(2.6,3)(6.8,3)
  \psline(4.5,0.6)(4.5,3.4)
  \psline(3,0.6)(3,3.4)
  \psline[linestyle=dashed](0.8,1)(3.4,1)
  \rput(7.5,-0.4){$1$}
  \rput(9.4,2){$0$}
  \rput(8.1,4.4){$(\frac{2}{3})$}
  \rput(7.2,3.5){$(\frac{1}{3})$}
  \rput(6.6,2.5){$-3$}
  \rput(3.75,3.4){$(-\frac{2}{3})$}
  \rput(5,2){$(-\frac{1}{3})$}
  \rput(3.5,2){$(-\frac{4}{3})$}
  \rput(2,1.4){$(-2)$}
  \rput(1.2,0.6){$C$}
  \rput(5,-1){VI $(n=6)$ (b)}
\end{pspicture*}\psset{xunit=1.15}
\begin{pspicture*}(-1,-1.4)(7,4.8)
  \psline(-0.4,0)(5.8,0)
  \psline(0,-0.4)(0,2.4)
  \psline(-0.4,2)(2.4,2)
  \psline[linestyle=dashed](5.4,-0.4)(5.4,2.4)
  \psline(2,1.6)(2,4.4)
  \psline(1.6,4)(4.4,4)
  \psline(3.75,4.25)(5.7,1.6)
  \psline(1,1.6)(1,4.4)
  \psline[linestyle=dashed](-1,3)(1.4,3)
  \psline[linestyle=dashed](4.6,2.2)(7,4)
  \rput(2.7,-0.4){$0$}
  \rput(-0.4,1){$(\frac{2}{3})$}
  \rput(0.5,1.6){$(\frac{1}{3})$}
  \rput(0.5,3.8){$(-\frac{1}{3})$}
  \rput(-0.5,3.4){$(-1)$}
  \rput(2.4,3){$(\frac{1}{3})$}
  \rput(4.8,3.5){$(\frac{1}{3})$}
  \rput(7.4,0.5){$(\frac{1}{3})$}
  \rput(3,4.4){$(\frac{1}{3})$}
  \rput(6.7,3.2){$C$}
  \rput(3,3.6){$-3$}
  \rput(6,3.8){$(0)$}
  \rput(3,-1){IX $(n=6)$}
\end{pspicture*}
\end{center}
\caption{$-f^*(K_{\bar X})$\quad $(c_1+c_2+r=0)$}\label{fig1.4}
\end{figure}
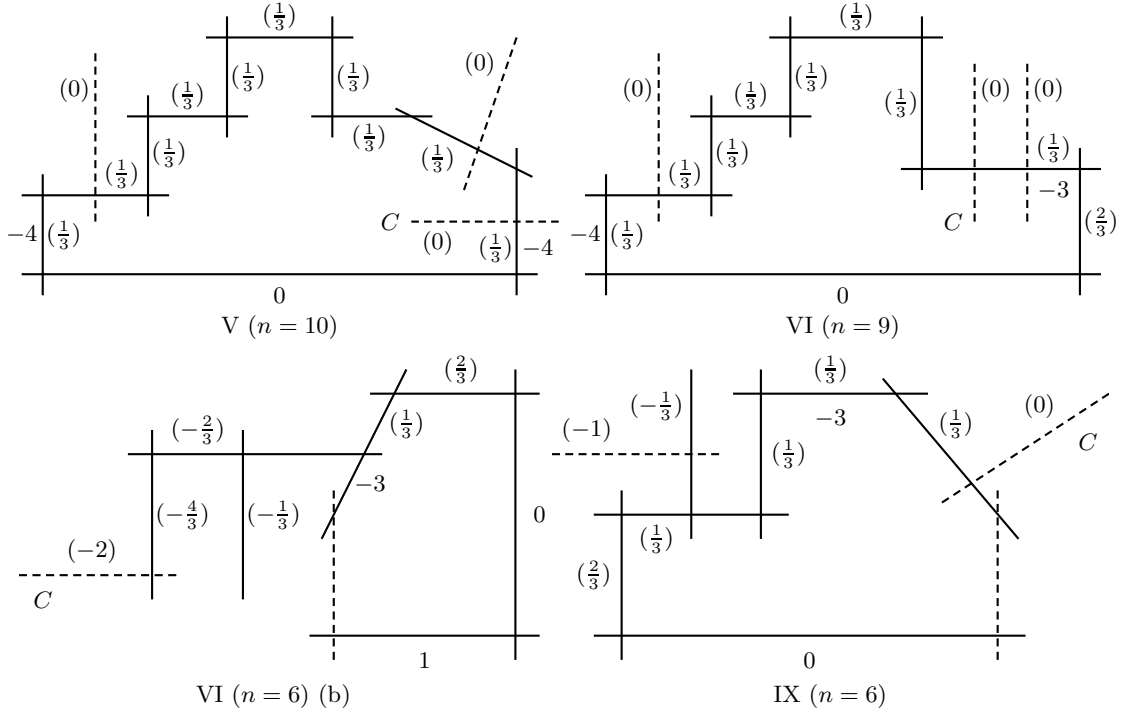

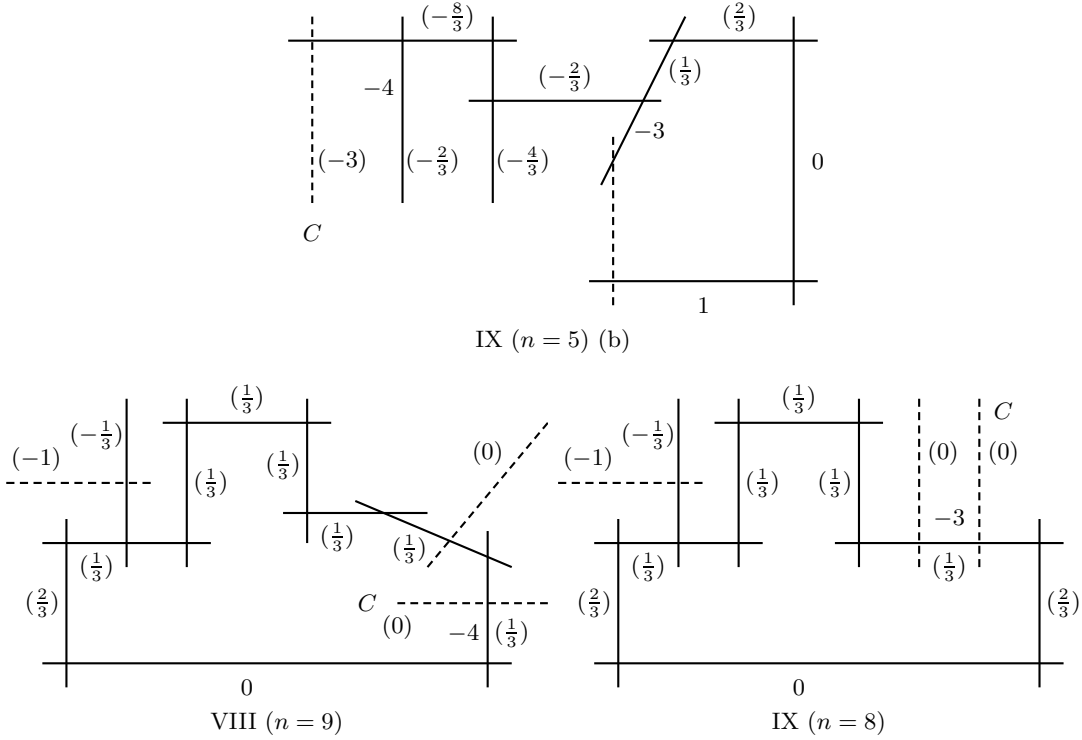
\begin{figure}[h!]
\begin{center}\psset{dash=3pt 2pt,unit=0.8}\everypsbox{\scriptsize}
\begin{pspicture*}(0,-1.4)(10,4.8)
  \psline(5.6,0)(9.4,0)
  \psline(9,-0.4)(9,4.4)
  \psline(6.6,4)(9.4,4)
  \psline[linestyle=dashed](6,-0.4)(6,2.4)
  \psline(5.8,1.6)(7.2,4.4)
  \psline(3.6,3)(6.8,3)
  \psline(4,1.3)(4,4.4)
  \psline[linestyle=dashed](1,1.3)(1,4.4)
  \psline(2.5,1.3)(2.5,4.4)
  \psline(0.6,4)(4.4,4)
  \rput(7.5,-0.4){$1$}
  \rput(9.4,2){$0$}
  \rput(8.1,4.4){$(\frac{2}{3})$}
  \rput(7.2,3.5){$(\frac{1}{3})$}
  \rput(6.6,2.5){$-3$}
  \rput(5.2,3.4){$(-\frac{2}{3})$}
  \rput(4.5,2){$(-\frac{4}{3})$}
  \rput(3.25,4.4){$(-\frac{8}{3})$}
  \rput(3,2){$(-\frac{2}{3})$}
  \rput(1.5,2){$(-3)$}
  \rput(2.1,3.2){$-4$}
  \rput(1,0.8){$C$}
  \rput(5,-1){IX $(n=5)$ (b)}
\end{pspicture*}
\end{center}
\begin{center}\psset{dash=3pt 2pt,unit=0.8}\everypsbox{\scriptsize}
\begin{pspicture*}(-1,-1.4)(8,4.8)
  \psline(-0.4,0)(7.4,0)
  \psline(0,-0.4)(0,2.4)
  \psline(-0.4,2)(2.4,2)
  \psline(3.6,2.5)(6,2.5)
  \psline(7,-0.4)(7,2.2)
  \psline(2,1.6)(2,4.4)
  \psline(4,2)(4,4.4)
  \psline(1.6,4)(4.4,4)
  \psline(1,1.6)(1,4.4)
  \psline(4.8,2.7)(7.4,1.6)
  \psline[linestyle=dashed](-1,3)(1.4,3)
  \psline[linestyle=dashed](5.5,1)(8,1)
  \psline[linestyle=dashed](6,1.6)(8,4)
  \rput(3,-0.4){$0$}
  \rput(-0.4,1){$(\frac{2}{3})$}
  \rput(0.5,1.6){$(\frac{1}{3})$}
  \rput(0.5,3.8){$(-\frac{1}{3})$}
  \rput(-0.5,3.4){$(-1)$}
  \rput(2.4,3){$(\frac{1}{3})$}
  \rput(3.6,3.3){$(\frac{1}{3})$}
  \rput(7.4,0.5){$(\frac{1}{3})$}
  \rput(3,4.4){$(\frac{1}{3})$}
  \rput(6.6,0.5){$-4$}
  \rput(5,1){$C$}
  \rput(5.5,0.6){$(0)$}
  \rput(5.7,1.9){$(\frac{1}{3})$}
  \rput(7,3.5){$(0)$}
  \rput(4.5,2.1){$(\frac{1}{3})$}
  \rput(3.5,-1){VIII $(n=9)$}
\end{pspicture*}
\begin{pspicture*}(-1,-1.4)(8,4.8)
  \psline(-0.4,0)(7.4,0)
  \psline(0,-0.4)(0,2.4)
  \psline(-0.4,2)(2.4,2)
  \psline(3.6,2)(7.4,2)
  \psline(7,-0.4)(7,2.4)
  \psline(2,1.6)(2,4.4)
  \psline(4,1.6)(4,4.4)
  \psline(1.6,4)(4.4,4)
  \psline(1,1.6)(1,4.4)
  \psline[linestyle=dashed](-1,3)(1.4,3)
  \psline[linestyle=dashed](5,1.6)(5,4.4)
  \psline[linestyle=dashed](6,1.6)(6,4.4)
  \rput(3,-0.4){$0$}
  \rput(-0.4,1){$(\frac{2}{3})$}
  \rput(0.5,1.6){$(\frac{1}{3})$}
  \rput(0.5,3.8){$(-\frac{1}{3})$}
  \rput(-0.5,3.4){$(-1)$}
  \rput(2.4,3){$(\frac{1}{3})$}
  \rput(3.6,3){$(\frac{1}{3})$}
  \rput(5.5,1.6){$(\frac{1}{3})$}
  \rput(5.5,2.4){$-3$}
  \rput(7.4,1){$(\frac{2}{3})$}
  \rput(3,4.4){$(\frac{1}{3})$}
  \rput(5.4,3.5){$(0)$}
  \rput(6.4,3.5){$(0)$}
  \rput(6.4,4.2){$C$}
  \rput(3.5,-1){IX $(n=8)$}
\end{pspicture*}
\end{center}
\caption{$-f^*(K_{\bar X})$\quad $(c_1+c_2+r<0)$}\label{fig1.5}
\end{figure}

The weighted dual graphs for some $-f^*(K_{\bar X})$ are illustrated
in Figures~\ref{fig1.4} and \ref{fig1.5}. For each of the irreducible
curves, the label with brackets indicates its coefficient, and the
label without brackets indicates its self-intersection number. The
labels for coefficient 1 are omitted.  A dotted line stands for a
$(-1)$-curve, and a solid line stands for a $(-2)$-curve if its
self-intersection number is not indicated.

\begin{Proof}[Proof of Lemma~\textup{\ref{l1.12}}]
 From the proof of Proposition~\ref{p1.9} (2),
\begin{equation*}
  (-K_{\bar X})^2= (K_X)^2+ \sum_{j=1}^n a_j(-2-(D_j)^2) = 8-n+ \sum_{j=1}^n a_j(-2-(D_j)^2).
\end{equation*}
The size $n$ of $D$ in Figure~\ref{fig1.1} is chosen so that
$n>8+\sum_{j=1}^n a_j(-2-(D_j)^2)$. Then $(-K_{\bar X})^2>0$. So by
Nakai-Moishezon criterion, $-K_{\bar X}$ is ample if and only
 $-K_{\bar X}\cdot \bar G>0$ for every irreducible curve $\bar G$ on $\bar X$.

 Let $\bar G$ be an irreducible curve on $\bar X$, and $G$ the proper transform of $\bar G$ on $X$. Then
\begin{equation*}
  -K_{\bar X}\cdot \bar G=-f^*(K_{\bar X})\cdot f^*(\bar G)= -f^*(K_{\bar X})\cdot G.
\end{equation*}
We will show that this number is positive by considering the
following two possibilities:

\medskip

\paragraph{\bf $G$ is contained in a fiber.}~

\emph{Case 1}. Suppose $G$ is a general fiber. Then $\bar G$ does
not contain the image of any center of blowup. So $G$ intersects
with the proper transforms of $M_r$ and $N_r$ on $X$. It follows
that $-f^*(K_{\bar X})\cdot G\geq 1+1/3>0$.

\emph{Case 2.} Suppose $G$ is contained in a singular fiber. Then
$G^2<0$. Note that $G\nsubseteq \Supp(D)$. By Lemma~\ref{l1.6}, $G$
is a $(-1)$-curve. Its coefficient in $-f^*(K_{\bar X})$ is the same
as that in $-K_X$.

(i) If $G\subseteq \Supp(\Delta_+)$, then $G$ intersects with
exactly two irreducible components of $\Delta$, which are contained
in $\Delta_+$. Moreover, exactly one of them is an irreducible
component of $D$. We have $-f^*(K_{\bar X})\cdot G\geq
(-1)+1/3+1>0$.

(ii) If $G\nsubseteq \Supp(\Delta_+)$, let $c$ be the coefficient of
$G$ in $-K_X$, then $G$ intersects with exactly one irreducible
component of $D$, whose coefficient in $-K_X$ is $c+1$. Note that
$G$ is disjoint from any other irreducible component of $\Delta$. So
$-f^*(K_{\bar X})\cdot G\geq (-1)c+(c+1-2/3)>0$.

\medskip

\paragraph{\bf $G$ is not contained in a fiber.}~

Note that $G_0:=\Phi(G)$ is a curve in $\mathbb F_r$. Write $G_0\sim
aM_r+bF_1$, where $a>0$ and $b\geq ar$. We have $G_0\cdot F_1=
G_0\cdot F_2=a$, $G_0\cdot M_r=b-ar\geq 0$ and $G_0\cdot N_r=b$. Let
$c_i$ be the smallest coefficient among all the irreducible
components of $\Phi^{-1}(F_i)$ appearing in $-f^*(K_{\bar X})$,
$i=1,2$. Then
\begin{equation}\label{eqn1.3}
  -f^*(K_{\bar X})\cdot G
  \geq ac_1+ac_2+0+b \geq a(c_1+c_2+r).
\end{equation}
By considering the sign of $c_1+c_2+r$, we have the following three
cases:

\emph{Case 1.} $c_1+c_2+r>0$. This is true for 22 configurations in
Figure~\ref{fig1.6}. For this case, it follows immediately from
(\ref{eqn1.3}) that $-f^*(K_{\bar X})\cdot G>0$.

\emph{Case 2.} $c_1+c_2+r=0$. There are 4 configurations as given in
Figure~\ref{fig1.4}.

For this case, we may assume that $b=ar$; otherwise $b>ar$ and
(\ref{eqn1.3}) implies that $-f^*(K_{\bar X})\cdot G\geq
a(c_1+c_2+r)+(b-ar)>0$. Then $G_0\sim aN_r$, and thus $G_0$ is
disjoint from the minimal section $M_r$. Therefore, there must exist
irreducible curves $L_i\subseteq \Phi^{-1}(F_i)$ with coefficient
$c_i$ appearing in $-f^*(K_{\bar X})$ such that $\Phi(L_i)$ is not a
point in $M_r$ $(i=1,2)$. However, it is easy to see from
Figure~\ref{fig1.4} that $F_1$ does not exist for any of these 4
configurations.

\emph{Case 3.} $c_1+c_2+r<0$. There are 3 configurations as given in
Figure~\ref{fig1.5}.

For each of them, denote $\{P_i\}:=M_r\cap F_i$ $(i=1,2)$, and let
$C', C''$ be the irreducible curves in $\Phi^{-1}(F_1)$ with
coefficients $\leq -(c_2+r)$ in $-f^*(K_{\bar X})$. Suppose that
$-f^*(K_{\bar X})\cdot G\leq 0$. Then $s:=(C'+C'')\cdot G>0$.

(i) VI $(n=6)$ (b). By computing the multiplicities of the center of
blowups, we have $(F_1\cdot G_0)_{P_1}\geq 4s$ and $(M_1\cdot
G_0)_{P_1}\geq 4s$. In particular, $G_0\sim aM_1+bF_1$ with $a\geq
4s$ and $b\geq 8s$. Then it would follow that $-f^*(K_{\bar X})\cdot
G\geq (-3)s+4s+8s>0$, a contradiction.

(ii) and (iii). VIII $(n=9)$ and IX $(n=8)$. For these cases,
$(M_0\cdot G_0)_{P_1}\geq s$ and $(F_1\cdot G_0)_{P_1}\geq 2s$. If
$P_2\in F_2\cap G_0$, then $G_0\cdot N_0\geq (G_0\cdot
M_0)_{P_1}+(G_0\cdot M_0)_{P_2}\geq s+1$. We would have
$-f^*(K_{\bar X})\cdot G\geq (-1)s+(s+1)>0$. Suppose $P_2\notin
F_2\cap G_0$.

IX $(n=8)$: Let $F_2'$ be the proper transform of $F_2$ on $X$. Then
$G\cdot F_2'= G_0\cdot F_2\geq 2s$. But then $-f^*(K_{\bar X})\cdot
G\geq (-1)s+(2/3)2s+s>0$, a contradiction.

VII $(n=9)$: Note that $-f^*(K_{\bar X})\cdot G\geq (-1)s+s=0$. If
$-f^*(K_{\bar X})\cdot G=0$, then $G_0\cdot M_0=(G_0\cdot
M_0)_{P_1}=s$ and $G_0\cdot F_1=(G_0\cdot F_1)_{P_1}=2s$; that is,
$G_0\sim 2sM_0+sF_1$. Note that $G$ is disjoint from $F_2'$. Then
$G\cdot C=2s-G\cdot F_2'=2s$. However, this would imply that $G_0$
has multiplicity $2s$ at the point $\Phi(C)$, and thus $s=G_0\cdot
M_0\geq 2s$, a contradiction again.

Therefore, $-K_{\bar X}\cdot \bar G=f^*(K_{\bar X})\cdot G>0$ for
every irreducible curve $\bar G$ on $\bar X$. Since $(-K_{\bar
X})^2>0$, by Nakai-Moishezon criterion, $-K_{\bar X}$ is ample for
all the 29 configurations listed in Figure~\ref{fig1.6}. We have
completed the proof of Lemma~\ref{l1.12}. \qed
\end{Proof}

\section[List of Weighted Dual Graphs]{The List of Weighted Dual Graphs of $C+D$}

{\footnotesize \psset{dotstyle=o,linewidth=0.5pt,unit=0.75}
\renewcommand{\baselinestretch}{1} \parskip=1mm

I (a) and (b):
\begin{center}\everypsbox{\scriptsize}
\begin{pspicture*}(-1,-0.5)(2,0.1)
  \psdots(0,0)(1,0)
  \uput[-90](0,0){$-3$}
  \uput[-90](1,0){$-1$}
\end{pspicture*}
\end{center}

\begin{center}\everypsbox{\scriptsize}
\begin{pspicture*}(-1,-0.5)(2,0.1)
  \psline[arrows=o-o](0,0)(1,0)
  \uput[-90](0,0){$-3$}
  \uput[-90](1,0){$-1$}
\end{pspicture*}
\end{center}

II (a) and (b):
\begin{center}\everypsbox{\scriptsize}
\begin{pspicture*}(-1,-0.5)(2,0.1)
  \psdots(0,0)(1,0)
  \uput[-90](0,0){$-6$}
  \uput[-90](1,0){$-1$}
\end{pspicture*}
\end{center}

\begin{center}\everypsbox{\scriptsize}
\begin{pspicture*}(-1,-0.5)(2,0.1)
  \psline[arrows=o-o](0,0)(1,0)
  \uput[-90](0,0){$-6$}
  \uput[-90](1,0){$-1$}
\end{pspicture*}
\end{center}

III:
\begin{center}\everypsbox{\scriptsize}
\begin{pspicture*}(-1,-0.5)(3,0.1)
  \psline[arrows=o-](0,0)(1,0)
  \psline[arrows=o-o](1,0)(2,0)
  \uput[-90](0,0){$-1$}
  \uput[-90](1,0){$-2$}
  \uput[-90](2,0){$-5$}
\end{pspicture*}
\end{center}

V $(n=5)$ (a) and (b):
\begin{center}\everypsbox{\scriptsize}
\begin{pspicture*}(-1,-0.5)(6,0.1)
  \psline[arrows=o-](0,0)(1,0)
  \psline[arrows=o-](1,0)(2,0)
  \psline[arrows=o-](2,0)(3,0)
  \psline[arrows=o-](3,0)(4,0)
  \psline[arrows=o-o](4,0)(5,0)
  \uput[-90](0,0){$-1$}
  \uput[-90](1,0){$-4$}
  \uput[-90](2,0){$-2$}
  \uput[-90](3,0){$-2$}
  \uput[-90](4,0){$-2$}
  \uput[-90](5,0){$-4$}
\end{pspicture*}
\end{center}

\begin{center}\everypsbox{\scriptsize}
\begin{pspicture*}(-1,-0.5)(5,1.2)
  \psline[arrows=o-](0,0)(1,0)
  \psline[arrows=-](1,0)(2,0)
  \psline[arrows=o-](2,0)(3,0)
  \psline[arrows=o-o](3,0)(4,0)
  \psline[arrows=o-o](1,0)(1,1)
  \uput[-90](0,0){$-4$}
  \uput[-90](1,0){$-2$}
  \uput[-90](2,0){$-2$}
  \uput[-90](3,0){$-2$}
  \uput[-90](4,0){$-4$}
  \uput[180](1,1){$-1$}
\end{pspicture*}
\end{center}

V $(n=6)$ (a) and (b):
\begin{center}\everypsbox{\scriptsize}
\begin{pspicture*}(-1,-0.5)(7,0.1)
  \psline[arrows=o-](0,0)(1,0)
  \psline[arrows=o-](1,0)(2,0)
  \psline[arrows=o-](2,0)(3,0)
  \psline[arrows=o-](3,0)(4,0)
  \psline[arrows=o-o](4,0)(5,0)
  \psdot(6,0)
  \uput[-90](0,0){$-4$}
  \uput[-90](1,0){$-2$}
  \uput[-90](2,0){$-2$}
  \uput[-90](3,0){$-2$}
  \uput[-90](4,0){$-2$}
  \uput[-90](5,0){$-4$}
  \uput[-90](6,0){$-1$}
\end{pspicture*}
\end{center}

\begin{center}\everypsbox{\scriptsize}
\begin{pspicture*}(-1,-0.5)(7,0.1)
  \psline[arrows=o-](0,0)(1,0)
  \psline[arrows=o-](1,0)(2,0)
  \psline[arrows=o-](2,0)(3,0)
  \psline[arrows=o-](3,0)(4,0)
  \psline[arrows=o-](4,0)(5,0)
  \psline[arrows=o-o](5,0)(6,0)
  \uput[-90](0,0){$-4$}
  \uput[-90](1,0){$-2$}
  \uput[-90](2,0){$-2$}
  \uput[-90](3,0){$-2$}
  \uput[-90](4,0){$-2$}
  \uput[-90](5,0){$-4$}
  \uput[-90](6,0){$-1$}
\end{pspicture*}
\end{center}

V $(n=10)$:
\begin{center}\everypsbox{\scriptsize}
\begin{pspicture*}(-1,-0.5)(11,0.1)
  \psline[arrows=o-](0,0)(1,0)
  \psline[arrows=o-](1,0)(2,0)
  \psline[arrows=o-](2,0)(3,0)
  \psline[arrows=o-](3,0)(4,0)
  \psline[arrows=o-](4,0)(5,0)
  \psline[arrows=o-](5,0)(6,0)
  \psline[arrows=o-](6,0)(7,0)
  \psline[arrows=o-](7,0)(8,0)
  \psline[arrows=o-](8,0)(9,0)
  \psline[arrows=o-o](9,0)(10,0)
  \uput[-90](0,0){$-4$}
  \uput[-90](1,0){$-2$}
  \uput[-90](2,0){$-2$}
  \uput[-90](3,0){$-2$}
  \uput[-90](4,0){$-2$}
  \uput[-90](5,0){$-2$}
  \uput[-90](6,0){$-2$}
  \uput[-90](7,0){$-2$}
  \uput[-90](8,0){$-2$}
  \uput[-90](9,0){$-4$}
  \uput[-90](10,0){$-1$}
\end{pspicture*}
\end{center}

VI $(n=4)$:
\begin{center}\everypsbox{\scriptsize}
\begin{pspicture*}(-1,-0.5)(4,1.2)
  \psline[arrows=o-](0,0)(1,0)
  \psline[arrows=-](1,0)(2,0)
  \psline[arrows=o-o](2,0)(3,0)
  \psline[arrows=o-o](1,0)(1,1)
  \uput[-90](0,0){$-2$}
  \uput[-90](1,0){$-3$}
  \uput[-90](2,0){$-2$}
  \uput[-90](3,0){$-4$}
  \uput[180](1,1){$-1$}
\end{pspicture*}
\end{center}

VI $(n=5)$ (a) and (b):
\begin{center}\everypsbox{\scriptsize}
\begin{pspicture*}(-1,-0.5)(6,0.1)
  \psline[arrows=o-](0,0)(1,0)
  \psline[arrows=o-](1,0)(2,0)
  \psline[arrows=o-](2,0)(3,0)
  \psline[arrows=o-](3,0)(4,0)
  \psline[arrows=o-o](4,0)(5,0)
  \uput[-90](0,0){$-1$}
  \uput[-90](1,0){$-2$}
  \uput[-90](2,0){$-3$}
  \uput[-90](3,0){$-2$}
  \uput[-90](4,0){$-2$}
  \uput[-90](5,0){$-4$}
\end{pspicture*}
\end{center}

\begin{center}\everypsbox{\scriptsize}
\begin{pspicture*}(-1,-0.5)(5,1.2)
  \psline[arrows=o-](0,0)(1,0)
  \psline[arrows=o-](1,0)(2,0)
  \psline[arrows=-](2,0)(3,0)
  \psline[arrows=o-o](3,0)(4,0)
  \psline[arrows=o-o](2,0)(2,1)
  \uput[-90](0,0){$-2$}
  \uput[-90](1,0){$-3$}
  \uput[-90](2,0){$-2$}
  \uput[-90](3,0){$-2$}
  \uput[-90](4,0){$-4$}
  \uput[180](2,1){$-1$}
\end{pspicture*}
\end{center}

VI $(n=6)$ (a) and (b):
\begin{center}\everypsbox{\scriptsize}
\begin{pspicture*}(-1,-0.5)(6,1.2)
  \psline[arrows=o-](0,0)(1,0)
  \psline[arrows=-](1,0)(2,0)
  \psline[arrows=o-](2,0)(3,0)
  \psline[arrows=o-](3,0)(4,0)
  \psline[arrows=o-o](4,0)(5,0)
  \psline[arrows=o-o](1,0)(1,1)
  \uput[-90](0,0){$-2$}
  \uput[-90](1,0){$-3$}
  \uput[-90](2,0){$-2$}
  \uput[-90](3,0){$-2$}
  \uput[-90](4,0){$-2$}
  \uput[-90](5,0){$-4$}
  \uput[180](1,1){$-1$}
\end{pspicture*}
\end{center}

\begin{center}\everypsbox{\scriptsize}
\begin{pspicture*}(-1,-0.5)(6,1.2)
  \psline[arrows=o-](0,0)(1,0)
  \psline[arrows=o-](1,0)(2,0)
  \psline[arrows=o-](2,0)(3,0)
  \psline[arrows=o-](3,0)(4,0)
  \psline[arrows=-o](4,0)(5,0)
  \psline[arrows=o-o](4,0)(4,1)
  \uput[-90](0,0){$-2$}
  \uput[-90](1,0){$-3$}
  \uput[-90](2,0){$-2$}
  \uput[-90](3,0){$-2$}
  \uput[-90](4,0){$-2$}
  \uput[-90](5,0){$-4$}
  \uput[180](4,1){$-1$}
\end{pspicture*}
\end{center}

VI $(n=7)$:
\begin{center}\everypsbox{\scriptsize}
\begin{pspicture*}(-1,-0.5)(8,0.1)
  \psline[arrows=o-](0,0)(1,0)
  \psline[arrows=o-](1,0)(2,0)
  \psline[arrows=o-](2,0)(3,0)
  \psline[arrows=o-](3,0)(4,0)
  \psline[arrows=o-](4,0)(5,0)
  \psline[arrows=o-](5,0)(6,0)
  \psline[arrows=o-o](6,0)(7,0)
  \uput[-90](0,0){$-1$}
  \uput[-90](1,0){$-2$}
  \uput[-90](2,0){$-3$}
  \uput[-90](3,0){$-2$}
  \uput[-90](4,0){$-2$}
  \uput[-90](5,0){$-2$}
  \uput[-90](6,0){$-2$}
  \uput[-90](7,0){$-4$}
\end{pspicture*}
\end{center}

VI $(n=9)$:
\begin{center}\everypsbox{\scriptsize}
\begin{pspicture*}(-1,-0.5)(9,1.2)
  \psline[arrows=o-](0,0)(1,0)
  \psline[arrows=-](1,0)(2,0)
  \psline[arrows=o-](2,0)(3,0)
  \psline[arrows=o-](3,0)(4,0)
  \psline[arrows=o-](4,0)(5,0)
  \psline[arrows=o-](5,0)(6,0)
  \psline[arrows=o-](6,0)(7,0)
  \psline[arrows=o-o](7,0)(8,0)
  \psline[arrows=o-o](1,0)(1,1)
  \uput[-90](0,0){$-2$}
  \uput[-90](1,0){$-3$}
  \uput[-90](2,0){$-2$}
  \uput[-90](3,0){$-2$}
  \uput[-90](4,0){$-2$}
  \uput[-90](5,0){$-2$}
  \uput[-90](6,0){$-2$}
  \uput[-90](7,0){$-2$}
  \uput[-90](8,0){$-4$}
  \uput[180](1,1){$-1$}
\end{pspicture*}
\end{center}

VII $(n=5)$ (a) and (b):
\begin{center}\everypsbox{\scriptsize}
\begin{pspicture*}(-1,-0.5)(5,1.2)
  \psline[arrows=o-](0,0)(1,0)
  \psline[arrows=-](1,0)(2,0)
  \psline[arrows=o-](2,0)(3,0)
  \psline[arrows=o-o](3,0)(4,0)
  \psline[arrows=o-o](1,0)(1,1)
  \uput[-90](0,0){$-2$}
  \uput[-90](1,0){$-3$}
  \uput[-90](2,0){$-2$}
  \uput[-90](3,0){$-3$}
  \uput[-90](4,0){$-2$}
  \uput[180](1,1){$-1$}
\end{pspicture*}
\end{center}

\begin{center}\everypsbox{\scriptsize}
\begin{pspicture*}(-1,-0.5)(5,1.2)
  \psline[arrows=o-](0,0)(1,0)
  \psline[arrows=-](1,0)(2,0)
  \psline[arrows=o-](2,0)(3,0)
  \psline[arrows=o-o](3,0)(4,0)
  \psline[arrows=o-o](2,0)(2,1)
  \uput[-90](0,0){$-2$}
  \uput[-90](1,0){$-3$}
  \uput[-90](2,0){$-2$}
  \uput[-90](3,0){$-3$}
  \uput[-90](4,0){$-2$}
  \uput[180](2,1){$-1$}
\end{pspicture*}
\end{center}

VII $(n=6)$ (a) and (b):
\begin{center}\everypsbox{\scriptsize}
\begin{pspicture*}(-1,-0.5)(7,0.1)
  \psline[arrows=o-](0,0)(1,0)
  \psline[arrows=o-](1,0)(2,0)
  \psline[arrows=o-](2,0)(3,0)
  \psline[arrows=o-](3,0)(4,0)
  \psline[arrows=o-o](4,0)(5,0)
  \psdot(6,0)
  \uput[-90](0,0){$-2$}
  \uput[-90](1,0){$-3$}
  \uput[-90](2,0){$-2$}
  \uput[-90](3,0){$-2$}
  \uput[-90](4,0){$-3$}
  \uput[-90](5,0){$-2$}
  \uput[-90](6,0){$-1$}
\end{pspicture*}
\end{center}

\begin{center}\everypsbox{\scriptsize}
\begin{pspicture*}(-1,-0.5)(7,0.1)
  \psline[arrows=o-](0,0)(1,0)
  \psline[arrows=o-](1,0)(2,0)
  \psline[arrows=o-](2,0)(3,0)
  \psline[arrows=o-](3,0)(4,0)
  \psline[arrows=o-](4,0)(5,0)
  \psline[arrows=o-o](5,0)(6,0)
  \uput[-90](0,0){$-1$}
  \uput[-90](1,0){$-2$}
  \uput[-90](2,0){$-3$}
  \uput[-90](3,0){$-2$}
  \uput[-90](4,0){$-2$}
  \uput[-90](5,0){$-3$}
  \uput[-90](6,0){$-2$}
\end{pspicture*}
\end{center}


%


VIII $(n=4)$:
\begin{center}\everypsbox{\scriptsize}
\begin{pspicture*}(-1,-0.5)(3,1.2)
  \psline[arrows=o-](0,0)(1,0)
  \psline[arrows=-](1,0)(2,0)
  \psline[arrows=o-o](1,0)(1,1)
  \psline[arrows=o-o](2,0)(2,1)
  \uput[-90](0,0){$-4$}
  \uput[-90](1,0){$-2$}
  \uput[-90](2,0){$-2$}
  \uput[180](1,1){$-2$}
  \uput[0](2,1){$-1$}
\end{pspicture*}
\end{center}

VIII $(n=5)$ (a) and (b):
\begin{center}\everypsbox{\scriptsize}
\begin{pspicture*}(-1,-0.5)(5,1.2)
  \psdot(0,0)
  \psline[arrows=o-](1,0)(2,0)
  \psline[arrows=o-](2,0)(3,0)
  \psline[arrows=-o](3,0)(4,0)
  \psline[arrows=o-o](3,0)(3,1)
  \uput[-90](0,0){$-1$}
  \uput[-90](1,0){$-4$}
  \uput[-90](2,0){$-2$}
  \uput[-90](3,0){$-2$}
  \uput[-90](4,0){$-2$}
  \uput[180](3,1){$-2$}
\end{pspicture*}
\end{center}

\begin{center}\everypsbox{\scriptsize}
\begin{pspicture*}(-1,-0.5)(5,1.2)
  \psline[arrows=o-o](0,0)(1,0)
  \psline[arrows=o-o](1,0)(2,0)
  \psline[arrows=o-](2,0)(3,0)
  \psline[arrows=-o](3,0)(4,0)
  \psline[arrows=o-o](3,0)(3,1)
  \uput[-90](0,0){$-1$}
  \uput[-90](1,0){$-4$}
  \uput[-90](2,0){$-2$}
  \uput[-90](3,0){$-2$}
  \uput[-90](4,0){$-2$}
  \uput[180](3,1){$-2$}
\end{pspicture*}
\end{center}


VIII $(n=9)$:
\begin{center}\everypsbox{\scriptsize}
\begin{pspicture*}(-1,-0.5)(9,1.2)
  \psline[arrows=o-](0,0)(1,0)
  \psline[arrows=o-](1,0)(2,0)
  \psline[arrows=o-](2,0)(3,0)
  \psline[arrows=o-](3,0)(4,0)
  \psline[arrows=o-](4,0)(5,0)
  \psline[arrows=o-](5,0)(6,0)
  \psline[arrows=o-](6,0)(7,0)
  \psline[arrows=-o](7,0)(8,0)
  \psline[arrows=o-o](7,0)(7,1)
  \uput[-90](0,0){$-1$}
  \uput[-90](1,0){$-4$}
  \uput[-90](2,0){$-2$}
  \uput[-90](3,0){$-2$}
  \uput[-90](4,0){$-2$}
  \uput[-90](5,0){$-2$}
  \uput[-90](6,0){$-2$}
  \uput[-90](7,0){$-2$}
  \uput[-90](8,0){$-2$}
  \uput[180](7,1){$-2$}
\end{pspicture*}
\end{center}

IX $(n=5)$ (a) and (b):
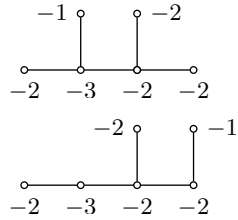
\begin{center}\everypsbox{\scriptsize}
\begin{pspicture*}(-1,-0.5)(4,1.2)
  \psline[arrows=o-](0,0)(1,0)
  \psline[arrows=-](1,0)(2,0)
  \psline[arrows=o-o](2,0)(3,0)
  \psline[arrows=o-o](1,0)(1,1)
  \psline[arrows=o-o](2,0)(2,1)
  \uput[-90](0,0){$-2$}
  \uput[-90](1,0){$-3$}
  \uput[-90](2,0){$-2$}
  \uput[-90](3,0){$-2$}
  \uput[180](1,1){$-1$}
  \uput[0](2,1){$-2$}
\end{pspicture*}
\end{center}

\begin{center}\everypsbox{\scriptsize}
\begin{pspicture*}(-1,-0.5)(4,1.2)
  \psline[arrows=o-](0,0)(1,0)
  \psline[arrows=o-](1,0)(2,0)
  \psline[arrows=-](2,0)(3,0)
  \psline[arrows=o-o](2,0)(2,1)
  \psline[arrows=o-o](3,0)(3,1)
  \uput[-90](0,0){$-2$}
  \uput[-90](1,0){$-3$}
  \uput[-90](2,0){$-2$}
  \uput[-90](3,0){$-2$}
  \uput[0](3,1){$-1$}
  \uput[180](2,1){$-2$}
\end{pspicture*}
\end{center}

IX $(n=6)$:
\begin{center}\everypsbox{\scriptsize}
\begin{pspicture*}(-1,-0.5)(6,1.2)
  \psline[arrows=o-](0,0)(1,0)
  \psline[arrows=o-](1,0)(2,0)
  \psline[arrows=o-](2,0)(3,0)
  \psline[arrows=o-](3,0)(4,0)
  \psline[arrows=-o](4,0)(5,0)
  \psline[arrows=o-o](4,0)(4,1)
  \uput[-90](0,0){$-1$}
  \uput[-90](1,0){$-2$}
  \uput[-90](2,0){$-3$}
  \uput[-90](3,0){$-2$}
  \uput[-90](4,0){$-2$}
  \uput[-90](5,0){$-2$}
  \uput[180](4,1){$-2$}
\end{pspicture*}
\end{center}

IX $(n=8)$:
\begin{center}\everypsbox{\scriptsize}
\begin{pspicture*}(-1,-0.5)(7,1.2)
  \psline[arrows=o-](0,0)(1,0)
  \psline[arrows=-](1,0)(2,0)
  \psline[arrows=o-](2,0)(3,0)
  \psline[arrows=o-](3,0)(4,0)
  \psline[arrows=o-](4,0)(5,0)
  \psline[arrows=-o](5,0)(6,0)
  \psline[arrows=o-o](1,0)(1,1)
  \psline[arrows=o-o](5,0)(5,1)
  \uput[-90](0,0){$-2$}
  \uput[-90](1,0){$-3$}
  \uput[-90](2,0){$-2$}
  \uput[-90](3,0){$-2$}
  \uput[-90](4,0){$-2$}
  \uput[-90](5,0){$-2$}
  \uput[-90](6,0){$-2$}
  \uput[180](1,1){$-1$}
  \uput[180](5,1){$-2$}
\end{pspicture*}
\end{center}
\begin{figure}[h!]
  \caption{Weighted Dual graphs of $C+D$}\label{fig1.6}
\end{figure}
}

\thanks
  \textsc{Acknowledgements}. The author would like to thank Prof D.-Q.\ Zhang for
introducing the topic and his kind guidance of the paper, and thank the referee for the valuable comments.

\bibliographystyle{amsplain}
\bibliography{Log}

\end{document}